
\input amstex
\documentstyle{amsppt}
\topmatter
\title
Cohomology, Symmetry, and Perfection
\endtitle
\thanks
This paper is dedicated to the memory of my friend Pere Menal.
I chose  the present topic for this occasion  because it was
the subject of our last conversation.  I miss him very much.\endgraf
To appear in the Fall 1992 issue, dedicated to the memory of Pere Menal,
of {\it Publicacions Matem\`atiques, } Universitat Aut\`onoma de
Barcelona\endthanks
\author Emili Bifet \endauthor
\address
Department of Mathematics
and Institute for Mathematical Sciences, State University of
New York, Stony Brook,
NY~~11794-3651
\endaddress
\email bifet\@ math.sunysb.edu \endemail
\endtopmatter
\document
\head 1.  Symmetry \endhead

In many situations that arise in  Algebraic Geometry one is
interested in  computing the multiplicative structure of the cohomology
ring $H^{\ast}(X)$  of some algebraic variety $X$.  Examples of such
situations include toric  varieties, complete quadrics, complete
symmetric varieties \cite{DP1-2}, \dots   Sometimes, as in
the examples just mentioned,
the variety $X$ is endowed with  symmetries  that reflect
the action of  some algebraic group $G$ on it.
In these cases there is a
recipe, inspired by the work of M\.~F\.~Atiyah and R\.~Bott \cite{AB1},
that
often works:

\roster
\item"1." Find a {\it strongly $G$-perfect}   decomposition
of $X$  (cf\.  Section $3$  below for a precise statement.)  In the
examples above, this is simply the  decomposition into
orbits. In general  there is a  natural candidate: the
Kempf-Hesselink  stratification of $X$ \cite{H,K,N}, a natural
outgrowth of D.~Mumford's Geometric Invariant Theory.

\item"2." With the help of this decomposition compute the equivariant
cohomology ring $H_G^*(X)$. (At this point it may also be
natural to apply the machinery of the localization theorem
\cite{AB2,Hs}; in so doing one usually obtains other interesting
descriptions of $H_G^*(X)$.)

\item"3." Recover $H^{\ast}(X)$ from $H_G^{\ast}(X)$.
\endroster

This recipe   is just one more instance
of the old philosophy of using any symmetries that may be
present in the problem in order to simplify it.

We shall work, for simplicity, with equivariant cohomology defined
in terms of homotopy quotients i.e. the Borel construction  (see
Section 2 below.)  But the knowledgeable reader could substitute
throughout $H_G^*(X)$ by  the $\ell$-adic cohomology of the algebraic
stack determined by  the $G$-variety $X$. He would thus gain the
advantage of having  a completely algebraic theory with, as a bonus,  a
very  interesting  arithmetic twist \cite{B2}.

In the examples mentioned earlier, the recipe works and gives
very explicit results. We shall  consider below in some detail the case
of   toric varieties, but the reader is advised to look at
\cite{BDP} for a thorough treatment of these  examples
from the present  point of view. In fact one of
the aims of this
paper is to better explain  the philosophy behind
those  computations   and to place them in a wider
conceptual setting.

Another aim of the paper is to outline  in the last section
an ``equivariant'' approach to some key results in the theory
of toric varieties.
This approach clarifies, I believe, the nature
of these results.

The text of the first three sections follows
closely a talk delivered
at the University of Copenhagen in
July 1989 on the occasion of
the Zeuthen Symposium. I should like to
thank S. Kleiman and A. Thorup
for organizing that conference and
creating a very friendly
athmosphere.

\head 2. Cohomology \endhead

Suppose  that a nice Lie group $G$ acts on a space $X$.  Na\"{\i}vely
speaking the  (equivariant)  cohomology of the $G$-space $X$
should be the
cohomology  of the quotient $X/G$. Unfortunately this gives
a useful
notion only  when $G$ is acting freely on $X$. In general it
is necessary
to find  the right notion of quotient, in fact a homotopy quotient
$X_G$,   in order to get a useful theory.
(~In a purely algebraic
context the role
of $X_G$ would be played by  the algebraic stack determined by the
$G$-variety $X$; the equivariant  cohomology  would just be the
$\ell$-adic cohomology of this stack.~) Before we can describe this
notion of quotient, however, it is necessary to look closely at the
case where $X$ is a point i.e. the theory of characteristic classes.

Recall that a {\it classifying space }
for principal
$G$-bundles is by definition a space $BG$ together with a
universal principal
$G$-bundle $EG$ over it.  By universal we mean that
isomorphism classes of
principal $G$-fibre bundles over a nice space $X$ correspond
naturally to
homotopy classes of maps from $X$ to $BG$.  The correspondence
is given by
pulling-back the universal bundle $EG$.  The cohomology
ring $H^{\ast}(BG)$
is by definition the {\it ring of characteristic classes } of $G$.

\example{Example 1}
$G = {\Bbb C}^{\times}$.  This is simply
the theory of line
bundles, and the classifying space is the infinite projective
space ${\Bbb
P}_{\Bbb C}^{\infty}$.  Its cohomology ring
$H^{\ast}(BG) = {\Bbb Z}[c_1]$ is a
polynomial ring in one variable of degree two.
\endexample
\example{Example 2}
$G = T = {\Bbb C}^{\times} \times \dots \times {\Bbb C}^{\times}$
(an algebraic torus.)  In this case we
have:
$$
BT \simeq B{\Bbb C}^{\times} \times \dots \times B{\Bbb C}^{\times}
$$
and its cohomology is a polynomial algebra in several variables
(as many as
factors.) In fact, if $X(T)$ denotes the group of algebraic characters
of $T$, then
$$
 H^{\ast}(BT) \simeq \text{Sym}^{\ast} X(T) \ .
$$
In particular $H^2(BT)\simeq X(T)$.
\endexample
\example{Example 3}
$G = \text{GL}_n({\Bbb C})$.  The classifying space is the infinite
dimensional  Grassmannian and
$$
H^{\ast}(BG) = {\Bbb Z}[c_1,\dots,c_n]
$$
where the variables $c_i$, $1 \le i \le n$, have degree $2i$ and
correspond to
the Chern classes.
\endexample

In general, if $T$ is a maximal torus in $G$, we have
$$
H^{\ast}(BG) = H^{\ast}(BT)^W
$$
where $W = N_G(T)/T$ is the Weyl group of $(G,T)$ and
the right hand side
denotes the subring of invariants.

We are now ready to describe the homotopy quotient    mentioned
earlier.
This  is given by
the {\it Borel construction} $X_G$ obtained after
exchanging the fibre $G$ of the universal bundle $EG$ with $X$ i.e.
$$
X_G = EG \times_G X = (EG \times X)/G
$$
where $G$ acts by $g \cdot (e,x) = (eg^{-1},gx)$.

\definition{Definition} The {\it equivariant cohomology }
$H_G^{\ast}(X)$ is by
definition the cohomology of the Borel construction $X_G$.
\enddefinition

Note that there is  a fibration
$$
X \to X_G \to BG. \tag 2.1
$$
The spectral sequence of this fibration  is the key to
the third step in the recipe. This is based on work of P. Deligne
\cite{De},  V. A. Ginzburg \cite{G}, F. Kirwan \cite{K}, \dots

Here follow some other properties  of equivariant cohomology
\cite{Hs,AB2}:

\roster
\item"a)" If $G$ acts freely on $X$, then
$H_G^{\ast}(X) = H^{\ast}(X/G)$.
\item"b)" If $T$ is a maximal torus in $G$
and $W = N_G(T)/T$ is the Weyl
group, then
$$
H_G(X) = H_T^{\ast}(X)^W.
$$
where the right hand side is the subring of $W$-invariants.
\item"c)" If $X$ has a single orbit, then
$$
H_G^{\ast}(X) \simeq H^{\ast}(BH).
$$
where $H$ is the stabilizer of any point.
\item"d)" If $K$ is a maximal compact subgroup of $G$, then
$$
H_G^{\ast}(X) \simeq H_K^{\ast}(X).
$$
\endroster

One of the reasons equivariant
cohomology is easier to compute than ordinary cohomology is that it
has many more ``points''. Let me try to explain this. Most
succesful calculations of cohomology achieve their objective by
expressing the cohomology of the space under consideration (e.g.
projective space $\Bbb{P}^n$) in terms of that of spaces for which
it is already known (e.g. cells.)  Ultimately, however, they reduce the
computation to that of the cohomology of a point.
If one thinks of
ordinary cohomology as being the case  $G=1$ of the equivariant one, then
it is clear that the points   coincide with the orbits.
Thus in the
equivariant theory every orbit gives rise to  a ``point'', and
there are as many points as there are conjugacy classes of subgroups in
$G$.  The equivariant cohomology of such a point $H$ is precisely the
ring of $H$-characteristic classes i.e. the cohomology of the
classifying space of $H$. It follows that in the equivariant theory
there is much more freedom of movement.

Another important feature of equivariant cohomology is that there
is a theory of equivariant Chern classes. A $G$-linearization of a
vector bundle $F$ over $X$ is an action $u:G\times F\to F$ which is linear
on the fibres and turns the projection $\pi :F\to X$ into a $G$-equivariant
map i.e. $\pi (g\cdot x)=g\cdot\pi (x)$ for
every $g\in G$ and every $x\in F$.
Note that the homotopy quotient $F_G$ provides us with a vector bundle
over $X_G$. The equivariant Chern classes of $(F,u)$ are by definition
the Chern classes of $F_G$.
This takes a most simple form for a line bundle over an orbit. In this
case the equivariant Chern class $c(L,u)$ is determined by the
isotropy action (character) of the stabilizer on the fibre of $L$ over
the point.
Actually these notions find their most natural formulation when expressed
in terms of algebraic stacks. For example a $G$-linearized
$\Cal O_X$-module is simply a module for the structure sheaf of the
algebraic stack determined by the $G$-variety $X$.

\head 3.  Perfection \endhead

Let $X$ be an algebraic variety, and let the algebraic group $G$
act on $X$.
Suppose $S \subset X$ is a closed $G$-invariant smooth subvariety
and let $U =
X - S$ be the complementary open set.  Under these conditions, there is
a long
exact sequence (the equivariant Thom--Gysin sequence, see \cite{AB1})
$$
\dots \to H_G^{i-2
\text{codim} S}(S) @>{\xi_S^i}>> H_G^i(X) @>>> H_G^i(U) \to
\dots~. \tag 3.1
$$
Moreover the composite of the maps
$$
\CD
H_G^{i-2\text{codim}S}(S) @>{\xi_S^i}>> H_G^i(X) \\
@.               @VV{\text{    restriction}}V \\
  {}    @.         H_G^i(S)
\endCD
$$
is multiplication by the Euler class $e(N_{S/X})$
(= top Chern class in this
context) of the normal bundle $N_{S/X}$.

If this long exact sequence splits into short exact sequences
$$
0 \to H_G^{i-2 \text{codim} S}(S) @>>> H_G^i(X)
@>>> H_G^i(U) \to 0 \tag 3.2
$$
then, for example, we have
$$
b_G^i(X) = b_G^i(U) + b_G^{i-2 \text{codim} S}(S)
$$
and one can deduce the equivariant Betti numbers of $X$ from
those of $S$ and
$U$.

In \cite{AB1} Atiyah and Bott made the following fundamental observation:

If $e(N_{S/X})$ is not a zero-divisor in the ring $H_G^*(S)$, then the
morphisms  $\xi_S^i$ are injective and the long exact sequences (3.1)
split  into short
exact sequences (3.2).

This motivates:

\definition{Definition}
We say that a decomposition
$$
X =   S_1\cup \ldots \cup S_N
$$
is {\it strongly $G$-perfect} if:
\roster
\item"1)" Each $S_i$ is both smooth and $G$-invariant.
\item"2)" For each $k$,
$$
X_k = S_1 \cup \dots \cup S_k
$$
is an open subset of $X$.
\item"3)" For each pair $(S_k,X_k)$, $k > 1$, the Euler class
$e(N_{S_k/X_k})$
is a non-zero divisor.
\endroster
\enddefinition

In \cite{AB1} a decomposition is defined to be
{\it $G$-perfect} if the
long exact sequence determined by each pair $(S_k,X_k)$ splits into
short exact sequences. In this case one has an identity of equivariant
Poincar\'{e}  series
$$
P^G_t(X)=\sum_{1\le i\le N} t^{2\cdot\text{codim}S_i}\
\cdot P^G_t(S_i) \
. $$
It is clear that   strongly $G$-perfect  implies   $G$-perfect.

An immediate consequence of the definition is

\proclaim{Proposition} If $\{ S_i \}_{1 \le i \le N}$ is a strongly
$G$-perfect decomposition of $X$, then for every partial union $X_k$
the morphism
induced  by the restrictions  to the strata
$$
H_G^{\ast}(X_k) @>>> \prod_{1 \le i \le k} H_G^{\ast}(S_i) \tag 3.3
$$
is injective.
\endproclaim

\demo{Proof} For $k = 1$, it is obvious.  Suppose it holds for $k - 1$;
it  suffices to show that the morphisms
$$
H_G^{\ast}(X_k) @>>> H_G^{\ast}(X_{k-1}) \times H^{\ast}(S_k)
$$
are injective.  Consider the diagram
$$
\CD
H_G^{i-2\text{codim} S_k}(S_k) @>{\xi}>> H_G^i(X_{k-1} \cup S_k) @>{\eta}>>
H_G^i(X_{k-1}) \\
@. @VV{\delta}V @. \\
{} @. H_G^i(S_k). @.
\endCD
$$
Now, if $\eta(a) = 0$, then there is a $b \in H_G^{i-2\text{codim}
S_k}(S_k)$  such that $a = \xi(b)$.  But, from
$$
0 = \d(a) = \d(\xi(b)) = b \cup e(N_{S_k/X_k})
$$
it follows that $b = 0$ and therefore $a = \xi(b) = 0$. \qed
\enddemo

Thus, in principle, if one knows the cohomology rings of the strata
and one controls the injection above, it is possible to describe the
cohomology ring of  $X$. This is the reason we singled out  this notion
for special consideration.

Atiyah and Bott  also give an infinitesimal criterion for $e(N_{S/X})$
to be a  non-zero divisor (see \cite{AB1} Proposition~13\.4.)
It is  proved in
\cite{K} using this criterion that the
Kempf-Hesselink stratification \cite{H}  of a $G$-variety is strongly
$G$-perfect in the  above sense.

Let us enunciate this last criterion  in the case of an orbit:

\proclaim{Proposition} Let ${\Cal O}$ be a $G$-orbit in the smooth
algebraic
variety $X$.  Choose a point $P \in {\Cal O}$ and
identify ${\Cal O}$ with
$G/G_P$, where $G_P$ is the stabilizer of $P$.  If the
isotropic action of a
maximal torus $T$ in $G_P$ on the normal space
$N_{{\Cal O}/X}(P)$ has no
non-zero fixed points, then $e(N_{{\Cal O}/X})$ is a
non-zero divisor.
\endproclaim

\demo{Proof} Since
$$
H_G^{\ast}({\Cal O}) = H_{G_P}^{\ast}(P) \hookrightarrow H_T^{\ast}(P)
$$
is an embedding, it suffices to see that $e = e(N_{{\Cal O}/X})$ is
non-zero
in $H_T^{\ast}(P)$.  But, if
$$
N_{{\Cal O}/X}(P) = \bigoplus_{\chi_i \in X(T)} \Bbb{C}_{\chi_i}
$$
is the weight decomposition, then
$$
e = \prod_i \chi_i \ne 0.
$$
since all $\chi_i$ are non-zero. \qed
\enddemo

The considerations above motivate:

\definition{Definition}
We say that $X$ is a {\it perfect embedding }
(regular embedding in
\cite{BDP}, but this was a bad choice to which I plead guilty) provided
\roster
\item"a)" Each orbit closure $\bar{\Cal O}$ is smooth and it
is the
transversal intersection of the codimension one orbit closures
that contain
it.
\item"b)" For every $P \in {\Cal O}$, the stabilizer $G_P$ has
a dense orbit
in the normal space $N_{{\Cal O}/X}(P)$.
\endroster
\enddefinition
To any perfect embedding $X$ we associate a simplicial
complex ${\Cal C}_X =
(V,S)$ as follows:
\roster
\item"1." $V = \{v \mid {\Cal O}_v$ is an orbit of codimension one\}.
\item"2." $\Gamma \subset V$ is a simplex if, and only if,
$$
\bigcap_{v \in \Gamma} \overline{{\Cal O}_v} \ne \varnothing.
$$
(Note that $\varnothing$ is a simplex.)
\endroster
It is easy to show that the simplexes are in one-to-one
correspondence with the orbits.

It is clear that the decomposition of $X$ into orbits is in this case
strongly $G$-perfect.
The algebraic  varieties mentioned    above, namely toric
varieties and complete
symmetric varieties (in particular complete quadrics),
provide examples of perfect embeddings.
In \cite{BDP} an explicit description, based on these ideas, is given
for the  equivariant cohomology ring   of any perfect embedding. It
should be possible to extend these  results to the case of a
well behaved
$G$-variety and the  Kempf-Hesselink stratification.

\head 4. An example: toric varieties \endhead

We shall illustrate the generalities of the preceding sections with
the  concrete  case of toric varieties.

Let $T$ be an algebraic torus. A {\it toric variety} is a normal
algebraic  $T$-variety $X$ containing $T$ as a dense orbit (this includes
the  requirement that the stabilizer at the points of this orbit be
trivial.) References \cite{D,F,O} provide very good expositions of the
basic theory  of these varieties.

A simple example is the affine plane with the action of
$T=\Bbb{C}^{\times}\times\Bbb{C}^{\times}$ given by
$$
(t_1,t_2)(x_1,x_2)=(t_1x_1,t_2x_2)\ .
$$
This action has four orbits, namely
the origin, the two punctured axes and the torus $T$ itself. This can be
generalized to the action
$$
t\cdot (x_1,x_2)=(t^{\chi_1}x_1,t^{\chi_2}x_2)
$$
determined by any integral basis $\{ \chi_1,\chi_2\}$ of the
character group $X(T)\simeq\Bbb{Z}^2$. This action induces an
action on
the ring of  regular functions on the affine plane, the ring of
polynomials in two  variables, given by $(t\cdot f)(x)=f(t^{-1}x)$. The
weight functions  (i.e. the functions $f$ such that $t\cdot f=t^{\chi}f$
for some  character $\chi$ called the weight  of $f$) are precisely the
monomials $\ a\cdot X_1^{n_1}X_2^{n_2}$ and their  weight is
$n_1(-\chi_1)+n_2(-\chi_2)$. These weights span a cone  $\sigma^{\vee}$
in
$X(T)\otimes \Bbb{R}$. Conversely we can recover the  variety, including
the $T$-action, by taking the Spectrum of  the monoid algebra $\Bbb{C}[
\sigma^{\vee}\cap X(T)]$ or, what is essentially the same, the
algebra homomorphisms from $\Bbb{C}[\sigma^{\vee}\cap X(T)]$ to
$\Bbb{C}$.

In general,
affine toric varieties can be constructed as follows.
We denote
$$
Y(T)=\text{Hom} (\Bbb{C}^{\times}, T)
$$
the group of one-parameter subgroups of the algebraic torus $T$.
(Recall that there is a pairing $X(T)\times Y(T)\to \Bbb Z$ given by taking
$\langle\chi,\mu\rangle$ to be the unique integer such that
$\chi (\mu (t)) =t^{\langle\chi,\mu\rangle}$
for all $t\in\Bbb C^{\times}$.)
First we
consider a  cone
$$
\sigma=\{ \lambda_1\mu_1 +\ldots +\lambda_m\mu_m \ |\
\lambda_i\in\Bbb{R},  \lambda_i\ge 0\text{ for every } i \}
$$
in $Y(T)_{\Bbb{R}}= Y(T)\otimes_{\Bbb{Z}}\Bbb{R}$ where $\{ \mu_1,\ldots,
\mu_m\}$ are finitely many one-parameter subgroups (actually we also ask
that the cone $\sigma$ have  a vertex i.e. $\sigma\cap(-\sigma)=0$.)
Next we introduce its dual in
$X(T)_{\Bbb{R}}=X(T)\otimes_{\Bbb{Z}}\Bbb{R}$ given by
$$
\sigma^{\vee}=\{ \chi\in X(T)_{\Bbb{R}}\  |\  \langle \chi ,\mu \rangle
\ge 0 \text{ for every } \mu\in\sigma\} \ .
$$
Then we express the monoid $\sigma^{\vee}\cap X(T)$ in terms of a
finite  number of generators
$$
\sigma^{\vee}\cap X(T)=\Bbb{N}\cdot{\chi_1 }+ \ldots +
\Bbb{N}\cdot{\chi_N} \tag 4.1
$$
(that this can always be done is a consequence of Gordan's lemma
\cite{D,F,O}.)
Finally we construct the affine model $X_{\sigma}$ by taking, as in the
case of the affine plane,  the Spectrum of  the monoid algebra or
equivalently the (scheme theoretic) closure of the image of
the map $T\to\Bbb{C}^N$ sending $t$ to
$(t^{\chi_1}, \ldots, t^{\chi_N})$.

A toric variety is obtained by glueing together the affine models
above along $T$-invariant open subsets. Think for example of the
projective plane obtained by glueing three copies of the affine plane
along the open orbit and identifying the punctured axes in pairs. Of
course, here  we are glueing not only the spaces but also the
$T$-actions so the actions on our three affine planes have to be
compatible. Fortunately the cones introduced earlier allow this
compatibility  to be expressed in simple terms.

\definition{Definition}
Let $T$ be an algebraic torus.
A  collection $\Sigma$ of cones as above in $Y(T)_{\Bbb R}$ is said to be a
{\it fan }  whenever it satisfies the following properties
\roster
\item"a)" Every face of a cone $\sigma$ in $\Sigma$  also belongs to  $\Sigma$.
\item"b)" The intersection of any two cones in $\Sigma$ is a face of
both of them.
\endroster
\enddefinition

A  fan $\Sigma$ gives rise to a toric variety $X_{\Sigma}$ obtained by
glueing together in turn  the pairs of  varieties $X_{\sigma_1}$,
$X_{\sigma_2}$ where
$\sigma_1, \sigma_2\in\Sigma$, along  the $T$-invariant open subset
$X_{\sigma_1\cap\sigma_2}$ determined by the common face
$\sigma_1\cap\sigma_2$.
 It is remarkable that there is a dictionary between  combinatorial
properties of  the fan   and geometric properties of $X_{\Sigma}$.
For example $X_{\Sigma}$ being compact (resp. smooth) translates into
the cones in $\Sigma$ covering all of $Y(T)_{\Bbb{R}}$ (resp. every
cone being such that  the $\chi_i$'s in (4.1) form a basis of
$X(T)$.) Moreover, there is a one-to-one correspondence between the
objects of the following three classes:
\roster
\item"1)" $T$-orbits.
\item"2)" $T$-invariant affine open subsets.
\item"3)" Cones in the fan.
\endroster
The correspondence is as follows. A cone $\sigma$ determines the
open subset $X_{\sigma}$, and this open set contains a unique
orbit $\Cal O_{\sigma}$  which is a closed subset in the relative
topology. Note that the codimension of an orbit coincides
with the dimension of the corresponding cone. In particular $T$-invariant
irreducible divisors $D_{\tau}=\overline {\Cal O_{\tau}}$ correspond
to one dimensional cones $\tau$.

It is easy to verify that any smooth toric variety $X_{\Sigma}$ gives a
perfect embedding and
we   associate to  it  a
simplicial complex
${\Cal C}_{\Sigma} = (V,S)$ defined as follows:

\roster
\item"i)" $V = \{\tau \in \Sigma: \dim \tau = 1\}$
\item"ii)" $S = \{\Gamma_{\sigma} \mid \sigma \in \Sigma\}$
where $\Gamma_{\sigma} = \{\tau \in V
\mid \tau \subset \sigma\}$.
\endroster
It is easy to see that
this simplicial complex is equivalent to the one
associated to $X_{\Sigma}$ as a perfect embedding.

Recall that the Reisner--Stanley algebra $R_{\Sigma}$ of ${\Cal C}_{\Sigma} =
(V,S)$ is by definition the quotient of the polynomial algebra
$$
{\Bbb Z}[x_v]_{v \in V}
$$
by the relations
$$
\prod_{v \in \Gamma} x_v = 0 \tag 4.2
$$
(one for each $\Gamma \notin S$.)

It is useful to think of the $x_v$'s as corresponding to the
$T$-invariant divisors $D_v$. Then the relation above reflects
the fact that $\bigcap_{v \in \Gamma} D_v =\varnothing$.
Similarly the monomials
$M=\prod x_v^{n_v}$ determine  a unique orbit closure
$\overline{\Cal O}=\bigcap_{n_v>0} D_v$
that we shall call the {\it support } of $M$.

Let $\Cal O_X(D_v)$ be the line bundle determined by the
effective divisor $D_v$. Recall that $D_v$ will coincide with the
zeros  of some section $s_v$ of $\Cal O_X(D_v)$.
We choose a $T$-linearization of each $\Cal O_X(D_v)$ in such a way
that in the induced action on sections $s_v$ has weight zero i.e.
it is $T$-invariant.

\proclaim{Proposition} Let $X_{\Sigma}$ be a smooth toric variety.
The natural morphism
$$
R_{\Sigma} \to H_T^{\ast}(X_{\Sigma})\ ,
$$
sending each $x_v$ to the equivariant Chern class of
$\Cal O_X(D_v)$,
is an
isomorphism.
\endproclaim
\demo{Proof}
I shall sketch the proof and refer the reader to
\cite{BDP}  for complete details. The idea is to define compatible
filtrations of both rings and to show that the induced homomorphism of
graded algebras is an isomorphism.

Choose an ordering of the orbits so that the partial unions
$X_k = \Cal O_1 \cup \dots \cup \Cal O_k$ are always open. Next we
filter
both $R_{\Sigma}$ and $H_T^*(X)$ as follows:
\roster
\item $F^kH_T^*(X)=\text{ker }H_T^*(X)\to H_T^*(X_{k-1})$.
\item $F^kR_{\Sigma}=
\text{ span of the monomials having support in }X-X_{k-1}$.
\endroster
The graded $k$-th pieces are respectively
$H_T^*(X_k,X_{k-1})\subset H_T^*(X_k)$
and the span $R_k$ of the monomials with support $\overline {\Cal O_k}$.
Note that $R_k=E_k\cdot \Bbb Z[x_v]_{D_v\supset\Cal O_k}$  where
$E_k=\prod_{D_v\supset{\Cal O_k}} x_v$  maps to the equivariant
Euler class of the normal bundle to $\overline{\Cal O_k}$. It
follows that
$R_k$ maps isomorphically onto the image of the Thom-Gysin map
$$
H^{*-2\text{codim}{\Cal O_k}}_T(\Cal O_k) @>>>
H_T^*(X_k,X_{k-1})\subset
H_T^*(X_k)
$$
and we are done. \qed
\enddemo

Another possible approach to this type of result is as follows. Every
affine open subset $X_{\sigma}$ retracts $T$-equivariantly to its closed
orbit $\Cal O_{\sigma}$. (This can be demonstrated very graphically in
the projective case using the compact torus
$T_{\text{compact}}=\text{U}(1)\times \dots\times\text{U}(1)$  and a  suitable
moment map.) Now a Mayer-Vietoris  argument allows us to conclude that
$$
\align
H_T^*(X)&= \varprojlim_{\sigma\in\Sigma} H_T^*(X_{\sigma})
= \varprojlim_{\sigma\in\Sigma} H_T^*(\Cal O_{\sigma}) = \\
        &= \varprojlim_{\sigma\in\Sigma} H^*( BT_{\sigma})
\endalign
$$
where we write  $T_{\sigma}$ for the stabilizer of any
point in  $\Cal O_{\sigma}$ and the
inverse systems on the right hand side
are indexed by the  inclusion relations among
the cones in the fan. (Note that $T_{\sigma}\subset T_{\eta}$ whenever
$\sigma\subset\eta$.)
These inverse limits describe  explicitly the image of
the restriction morphism  (3.3).

It is easy to see that the
extension of the morphism
$X(T) \to R_{\Sigma}$ given by
$$
\chi \mapsto \sum \Sb v \in V \\ v = (\mu_v) \endSb \langle\chi,\mu_v\rangle
\cdot
x_v  \tag 4.3
$$
can be identified with the natural morphism
$$
H^{\ast}(BT) \to H_T^{\ast}(X_{\Sigma}).
$$
When $X_{\Sigma}$ is  projective, the spectral
sequence determined by the
fibration (2.1) collapses (one argument:
the cohomology of base and fibre
vanish in  odd degrees \cite{BB1-3}.)  It follows that
we can identify the ring $H^*(X_{\Sigma})$
with the quotient of
$H_T^*(X_{\Sigma})$ by the ideal generated by
the strickly  positive part
of $H^*(BT)$. Similar  results hold  for
$X_{\Sigma}$ compact and
cohomology with  rational coefficients.
Thus one reobtains as was to be
expected the results of Jurkiewicz and Danilov \cite{D}.
But I believe
the present
approach clarifies the nature of the relations that appear
in that description.

Let me describe  from this point of view  the  equivariant
approach to the Picard group of $X_{\Sigma}$ that
I discovered while working on  \cite{B1, Remarks 2.1}.
The Picard group of $X$ can be identified with $H^2(X)$ \cite{BB1-3},
and it is clear
from the results above that since
$H^2(BT)=X(T)$,  we have
$$
H^2(X)=H^2_T(X)/X(T)
$$
where
$$
\align
H^2_T(X) &=\varprojlim_{\sigma\in\Sigma} H^2_T(\Cal O_{\sigma})  \\
         &=\varprojlim_{\sigma\in\Sigma} X(T_{\sigma})  \ .
\endalign
$$
But taking characters in the short exact sequence
$$
1\to  T_{\sigma}\to T\to T/T_{\sigma}\to 1
$$
gives
$$
X(T_{\sigma})= X(T)/(\sigma^{\perp}\cap X(T))
$$
as was to be expected. Actually, all this can be done purely
algebraically  using the equivariant Picard group $\text{Pic}_T(X)$.
This is the group of isomorphism classes of pairs $(L,u)$ consisting
of a line bundle $L$ over $X$  and a $T$-linearization of it.
The usual theory of line bundles on toric
varieties can now be realized as follows. Every line bundle $L$ on $X$
has a $T$-linearization and its equivariant Chern class $c(L,u)$ is given
by the family
$(\chi_{\sigma})_{\sigma\in\Sigma}\in H^2_T(X)$
where
each $\chi_{\sigma}\in X(T_{\sigma})$ denotes the character describing
the isotropy action of the stabilizer $T_{\sigma}$ on the fibre of $L$
over a  point in the orbit $\Cal O_{\sigma}$.
(This gives an isomorphism of $\text{Pic}_T(X)$ and $H^2_T(X)$.)
The cohomology spaces
$H^i(X,L)$ can now be completely described, as $T$-representations,
in terms of the equivariant
Chern class $c(L,u)$.
This is part of a more general phenomenon that includes the yoga of
moment maps. (The knowledgeable reader will recall
that the moment map is really a representative for the
equivariant Chern class
in the de~Rham model of equivariant cohomology \cite{AB2}.)
But this will have to be the object of another paper\dots


\enddocument